\documentclass{amsart}
\usepackage{amsmath}
\newtheorem{theorem}{Theorem}
\newtheorem{cor}[theorem]{Corollary}
\theoremstyle{definition} 
\newtheorem{example}[theorem]{Example}
\newtheorem*{remark}{Remark}

\newcommand{\fF}{{\mathfrak F}} 
\newcommand{\cser}{\mathcal{C}}
\newcommand{\Lie}{{}^{\text{Lie}} \mspace{-1.5mu}{}}

\title{Chief factors covered by projectors of soluble Leibniz algebras}
\author{Donald W. Barnes}
\address{1 Little Wonga Rd.\\Cremorne NSW 2090\\Australia\\}
\email{donwb@iprimus.com.au}

\subjclass[2010]{Primary 17A32, 17B30, Secondary 20D10}
\keywords{Leibniz algebras, Lie algebras, saturated formations, projectors}

\begin{document}

\begin{abstract} Let $\fF$ be a saturated formation of soluble Leibniz algebras.  Let $K$ be an $\fF$-projector and $A/B$ a chief factor of the soluble Leibniz algebra $L$.  It is well-known that if $A/B$ is $\fF$-central, then $K$ covers $A/B$.  I prove the converse: if $K$ covers $A/B$, then $A/B$ is $\fF$-central.
\end{abstract}

\maketitle

The theory of saturated formations and projectors for finite-dimensional soluble Lie and Leibniz algebras has been developed analogous to that for finite soluble groups.  In broad outline, the theories run parallel, but some results require very different proofs.  Not all the results translate.  One of the most troubling of these is the failure of the conjugacy theorem for projectors to translate to Lie and Leibniz algebras.  

A well-known basic theorem in all three contexts is that if $\fF$ is a saturated formation, $A/B$ an $\fF$-central chief factor of the group or algebra $L$ and if $K$ is an $\fF$-projector of $L$, then $K$ covers $A/B$, that is, $K+B \supseteq A$.  In this note, I prove the converse for Lie and Leibniz algebras.  For groups, if one $\fF$-projector covers the chief factor $A/B$ then trivially, every $\fF$-projector covers $A/B$ because of the conjugacy theorem.  Curiously, this result follows for Lie and Leibniz algebras as a consequence of the main theorem of this note,  a result which has no group theory analogue. 

All algebras considered in this note are soluble finite-dimensional Leibniz algebras over the field $F$ and $\fF$ is a saturated formation.  I denote the $\fF$-residual of the algebra $L$ by $L_{\fF}$.  If $A$ is an ideal of $L$, then $\cser_L(A)$ denotes the centraliser of $A$ in $L$.

\begin{theorem} \label{cover} Let $\fF$ be a saturated formation of soluble Leibniz algebras.  Let $A/B$ be a chief factor and $K$ an $\fF$-projector of the soluble Lie algebra $L$.  Suppose that $K$ covers $A/B$.  Then $A/B$ is $\fF$-central.
\end{theorem}

\begin{remark}   In view of  Barnes \cite[Corollary 3.17]{SchunckLeib}, the corresponding result for saturated formations of soluble Lie algebras may be regarded as a special case, one to which, in the proof, I reduce the general result.
\end{remark}

\begin{proof}  We can work in $L/B$, so we can assume that $B=0$ and $K \supseteq A$.  Put $L_1 = L/A$ and $K_1 = K/A$.  Then $K_1$ is a $\fF$-projector of $L_1$ and $A$ is an irreducible $L_1$-module which is $\fF$-hypercentral as $K_1$-module.  The assertion is equivalent to the assertion that $A$ is $\fF$-central as $L_1$-module.  As we can work with $L_1/\cser_{L_1}(A)$, we can assume that $\cser_{L_1}(A) = 0$, and so, that $\cser_L(A) = A$.  Thus $L$ is primitive and $A$ is the only minimal ideal of $L$.  By a theorem implicit in Loday and Pirashvili \cite{LodP} (see also Barnes \cite[Theorem 1.4]{EngelLeib}), $L_1$ is a Lie algebra and $A$ is either symmetric or antisymmetric.  By Barnes \cite[Theorem 3.16]{SchunckLeib}, we may suppose that $A$ is symmetric and so, that $L$ is a Lie algebra.  Further, an irreducible symmetric $L$-module is $\fF$-central if and only if it is $\Lie\fF$-central, where $\Lie\fF$ is the saturated formation of soluble Lie algebras consisting of those Lie algebras which are in $\fF$.

If $L/A \in \fF$, then $K = L$ and the result holds.  So we suppose that $L_{\fF}$ properly contains $A$.  From this, we shall derive a contradiction.

There exists an ideal $C$ such that $L_{\fF}/C$ is a chief factor of $L$.  Since $L/L_{\fF} \in \fF$ but $L/C \notin \fF$, the chief factor $L_{\fF}/C$ is an irreducible $\fF$-eccentric $K$-module.  Now consider $V = L_{\fF}/A$ as a $K$-module.   By Barnes \cite[Theorem 4.4]{HyperC}, we have a $K$-module direct decomposition $V = V^+ \oplus V^-$ where $V^+$ is $\fF$-hypercentral and $V^-$ is $\fF$-hypereccentric.  Since $L_{\fF}/C$ is $\fF$-eccentric, $V^- \ne 0$.

Consider the $K$-module homomorphism $\phi: V^- \otimes A \to A$ given by $\phi(\bar{v} \otimes a) = va$ for $\bar{v} = v+A \in V^-$ and $a \in A$.  Put $W = \phi(V^- \otimes  A)$.  Since $\cser_L(A) = A$, $W \ne 0$.  Since as $K$-modules, $V^-$ is $\fF$-hypereccentric and $A$ is $\fF$-hypercentral,  $V^- \otimes A$ is $\fF$-hypereccentric by Barnes \cite[Theorem 2.3]{hexc}.  Therefore $W$ is $\fF$-hypereccentric.  But $W$ is a submodule of the $\fF$-hypercentral module $A$. This contradiction completes the proof.
\end{proof}

\begin{cor} Suppose that some $\fF$-projector of $L$ covers the chief factor $A/B$.  Then every $\fF$-projector covers $A/B$.
\end{cor}

The group theory analogue of Theorem \ref{cover} is false.

\begin{example}  Let $\fF$ be the saturated formation of 2-groups.  It is locally defined by the family of formations $f(2)=\{1\}$, $f(p)= \emptyset$ for $p \ne 2$.  Then the $\fF$-projectors of a group $G$ are its $2$-Sylow subgroups.  Let $G = S_4$ be the symmetric group on four symbols and let $A$ be its normal subgroup of order $4$.  Then $A$ is covered by the $2$-Sylow subgroups of $G$ but is not $\fF$-central as $G/\cser_G(A) \notin f(2)$.
\end{example}

\bibliographystyle{amsplain}

\end{document}